\documentclass[11pt]{article}

\usepackage{amsmath,amssymb}
\usepackage{dsfont}
\usepackage{graphicx}
\usepackage{calrsfs}   
\usepackage{color,epsfig}
\definecolor{marin}{rgb}   {0.,   0.3,   0.7}
\definecolor{rouge}{rgb}   {0.8,   0.,   0.}
\definecolor{sepia}{rgb}   {0.8,   0.5,   0.}
\usepackage[colorlinks,citecolor=marin,linkcolor=rouge,
            bookmarksopen,
            bookmarksnumbered
           ]{hyperref}

\newtheorem{lemma}{Lemma}[section]

\newtheorem{proposition}[lemma]{Proposition}
\newtheorem{corollary}[lemma]{Corollary}
\newtheorem{remark}[lemma]{Remark}
\newtheorem{example}[lemma]{Example}

\newtheorem{notation}[lemma]{Notation}
\newtheorem{definition}[lemma]{Definition}
\newtheorem{conclusion}[lemma]{Conclusion}
\numberwithin{equation}{section}
\newcommand{\QED}{\mbox{}\hfill \raisebox{-0.2pt}{\rule{5.6pt}{6pt}\rule{0pt}{0pt}}
          \medskip\par}
\newenvironment{Proof}{\noindent    \parindent=0pt\abovedisplayskip = 0.5\abovedisplayskip
    \belowdisplayskip=\abovedisplayskip{\bfseries Proof. }}{\QED}

\newcommand{\dd}{\mathrm{d}}

\newcommand{\E}{\mathbb{E}}

\newcommand{\R}{\mathbb{R}}
\newcommand{\RR}{\rangle}
\newcommand{\LL}{\langle}
\newcommand{\Lr}{\mathrm{L}}

\newcommand{\Norm}[2]{\|#1\|\left.\vphantom{T_{j_0}^0}\!\!\right._{#2}}

\renewcommand{\P}{\mathbb{P}} 
 
\title{Reconciling alternate methods for the determination of charge distributions: A probabilistic approach
  to high-dimensional least-squares approximations.}

\begin{document}

\addtolength{\marginparwidth}{-15pt}

\author{Nicolas Champagnat$^1$, Christophe Chipot$^2$ and Erwan Faou$^3$}

\footnotetext[1]{TOSCA project-team, INRIA Sophia Antipolis --
  M\'editerran\'ee, 2004 route des Lucioles, BP.93, 06902 Sophia
  Antipolis Cedex, France (\texttt{Nicolas.Champagnat@sophia.inria.fr})}

\footnotetext[2]{Laboratoire de Chimie Th\'eorique, Unit\'e Mixte de
  Recherche CNRS no. 7565, Universit\'e Henri-Poincar\'e -- Nancy I,
  BP 239, 54506 Vandoeuvre-l\`es-Nancy Cedex, France
  (\texttt{Christophe.Chipot@edam.uhp-nancy.fr})}

\footnotetext[3]{INRIA \& Ecole Normale Sup\'erieure de Cachan
  Bretagne, Avenue Robert Schumann, 35170 Bruz, France
  (\texttt{Erwan.Faou@inria.fr})}

\maketitle

\abstract{ {We propose extensions and improvements} of the statistical analysis of distributed multipoles
  (SADM) {algorithm put forth} by {\sc Chipot} {\em et al} in \cite{chipot-98} for the derivation of
  {distributed} atomic multipoles from the {quantum-mechanical} electrostatic potential. {The method is
    mathematically extended to general least-squares problems and provides an alternative approximation method
    in cases where the original least-squares problem is computationally not tractable, either because of its
    ill-posedness or its high-dimensionality.} The solution is approximated {employing} a Monte Carlo method
  that takes the average of a random variable defined as the solutions of random small least-squares problems
  drawn as {subsystems} of the original problem. {The conditions that ensure convergence and consistency} of
  the method are {discussed}, {along with an analysis} of the computational cost in specific {instances}.


}


\section{Introduction}

In the realm of the molecular modeling of
complex chemical systems, atom-centered multipole distributions constitute a
popular route to simplify the description of intricate electron densities.
Streamlined down to their most rudimentary representation, these densities are
generally mimicked in macromolecular force fields by simple point charges, from
which, in the context of molecular simulations, Coulomb interactions can be
rapidly evaluated. Whereas nuclear charges are clearly centered onto the
constituent atoms, the electron charge distribution extends over the entire
molecular system. As a result, in sharp contrast with the higher-order
multipole moments of a neutral molecule, which, strictly speaking, are
quantum-mechanical observables, atomic charges cannot be defined univocally, in
an equally rigorous fashion. They ought to be viewed instead as a convenient
construct, the purpose of which is to reduce the complexity of molecular charge
distributions by means of compact sets of parameters providing a useful, albeit
naive framework to localize specific interactions onto atomic sites.

The ambiguous nature of atom-centered charges
has, therefore, prompted the development of alternative paths towards their
determination~\cite{Cornell1998}. The choice of the numerical scheme ought to
be dictated by three prevalent criteria, namely (i) the computational cost of
the derivation, (ii) the ease of implementation within the framework of a
physical model and (iii) the ability of the point-charge model to reproduce
properties of interest with the desired accuracy. Under a number of
circumstances, crude atomic charges determined through inexpensive calculations
are shown to be adequate. In other, more common scenarios, for instance, in
molecular simulations of complex chemical systems, the accurate description of
the electrostatic interactions at play can be of paramount importance. The
atomic charges utilized in these simulations are by and large derived from
quantum-mechanical calculations carried out at a reasonably high level of
theory, which in many cases, can be appreciably expensive. In the vast majority
of popular potential energy functions, point-charge models are derived
quantum-mechanically, following, in a nutshell, two distinct philosophies. On
the one hand, the numerical simulations of condensed phases imposes that
solute-solvent interactions be described as accurately as possible.
Accordingly, in macromolecular force fields like {\sc
Charmm}~\cite{MacKerell1998}, the atomic charges are determined based on a
series of independent quantum-mechanical calculations featuring different
relative positions of a solvent molecule around the solute. On the other hand,
the electrostatic potential can be viewed as the fingerprint of the molecule,
the accurate representation of which guarantees a reliable description of
intermolecular interactions. In potential energy functions like {\sc
Amber}~\cite{Cornell1995}, point charges are derived from the molecular
electrostatic potential, exploiting the fact that the latter is a
quantum-mechanical observable readily accessible from the wave function.

In their seminal article, Cox and
Williams~\cite{Cox1981} proposed an attractive approach, whereby sets of
atom-centered charges can be easily derived on the basis of a single-point
quantum-mechanical calculation. The electrostatic potential is evaluated on a
grid of $M$ points lying around the molecule of interest, outside the van der
Waals envelope of the latter. Restricting the multipole expansion of the
electrostatic potential to the monopole term, the charges borne by the $n$
atomic sites of the molecule are determined by minimizing the root-mean square
deviation between the reference, quantum-mechanical quantity and its
zeroth-order approximation --- i.e. $q_i T_{ki}^{00}$, where $T_{ki}^{00} =
\|x_i - x_k \|^{-1}$, is the potential created at point $k$ by atomic
site $i$. In its pioneering form, the algorithm handled the least-squares
problem iteratively. Chirlian and Francl subsequently proposed to resort to a
non-iterative numerical scheme~\cite{Chirlian1987}, which obviates the need for
initial guesses and solves the overdetermined system of linear equations
through matrix inversion. This route for the derivation of point-charge models
can be generalized in a straightforward fashion to higher-order multipoles.

The success of potential-derived charges stems
in large measure from their ease of computation and the demonstration for a
host of chemical systems that they are able to reproduce with an appreciable
accuracy a variety of physical properties. This success is, however, partially
clouded by one noteworthy shortcoming of the method --- point charges borne by
atoms buried in the molecule cannot be determined unambiguously from a
rudimentary least-squares fitting procedure. Symptomatically, for those
molecular systems, in which the contribution of the subset of buried atoms to
the electrostatic potential is ill-defined, the derived charges are in apparent
violation with the commonly accepted rules of electronegativity differences,
e.g. a C$^{\delta -}$---Cl$^{\delta +}$ bond polarity in carbon tetrachloride,
in lieu of the intuitive C$^{\delta +}$---Cl$^{\delta -}$. Bayly et al. tackled
this issue through the introduction of hyperbolic penalty functions in their
fitting procedure~\cite{Bayly1993}. Arguably enough, this numerical scheme
addresses the symptom rather than its actual cause. As was commented on by
Francl et al. in the light of singular-value-decomposition
analyses~\cite{Francl1996}, the matrices  of the least-squares problem are rank
deficient, to the extent that statistically valid charges cannot be assigned
univocally to the selected set of atoms in the molecule.

To delve further into this issue, Chipot {\em
et al.} proposed an alternative algorithm coined statistical analysis of
distributed multipoles (SADM)~\cite{chipot-98}, wherein atom-centered
multipoles are also derived from the quantum-mechanical electrostatic
potential, yet following a somewhat different pathway than the conventional
least-squares scheme. Instead of solving directly the $n \times M$
overdetermined system of linear equations, for instance through matrix
inversion, a subset of $n$ points is drawn amongst the $M$ points of the grid
and the corresponding $n \times n$ system of linear equations is solved. This
procedure, referred to as an \emph{experiment}, is repeated with different
subsets of grid points, from whence probability distributions are obtained for
the series of multipoles being sought. Strictly speaking, each probability
distribution ought to be determined from C$_M^n$ independent experiments. On
account of the computational burden, however
--- viz. typically, for a molecule formed by ten atoms and a grid of 2,000
points sampling the three-dimensional space around it, this would imply solving
approximately 2.76 $\times$ 10$^{26}$ systems of linear equations --- only 3--5
$\times$ 10$^5$ independent experiments are performed, which has proven
heuristically to be appropriate.


The mathematical description of this problem is the following: denoting by $(q_j)_{j = 1}^n$ the unknown
charges borne by the $n$ particles, and by $\gamma_j(x) = \|x- x_j\|^{-1}$, the electrostatic potential
associated with each $x_j \in \R^3$, the least square problem consists in finding the minimum $(q_j)_{j = 1}^n
\in \R^n$ of the function
\begin{equation}
\label{Eurpb}
\R^n \ni a \mapsto
\sum_{i = 1}^M | f(y_i) - \sum_{j = 1}^n a_j\gamma_j(y_i) |^2,
\end{equation}
where $(y_j)_{j = 1}^M \in \R^{3M}$ are the coordinates of the $M$ external points. Here $f(y_j)$ stands for
the approximation of the electrostatic potential {at $y_j$} obtained \color{black} by quantum-mechanical
{calculations.}

Instead of solving directly the problem \eqref{Eurpb}, the SADM consists in drawing $n$ points $ y^{(i)}$
amongst the $M$ points $y_j$, solve the $n \times n$ problem $f(y^{(i)}) = \sum_{j = 1}^n\gamma_j(y^{(i)})
a_j$, $i = 1,\ldots,n$ in the least squares sense and subsequently plot the distribution of each
$a_j$. In~\cite{chipot-98}, Chipot {\em et al.} notice that the {latter} are Cauchy-like distributions (with
seemingly infinite expectation) centered around the exact solution of the original least-squares problem.
Note that this method not only provides a numerical approximation of the solution, but {also} a global
statistical distribution {that reflects} the accuracy of the physical model {being utilized}.

Interestingly enough, it turns out that this kind of approach can be extended to many situations arising in
computational mathematics and physics. The principle of the SADM algorithm is in fact very general, and can be
adapted to derive efficient algorithms that are robust with the dimension of the underlying space of
approximation. This in turn provides new numerical methods of practical interest for high dimensional
{approximation} problems, where traditional {least-squares} methods are impossible to implement,
{either because of the high dimensionality or the ill-posedness of the least-squares problem.}

The goal of the present contribution is twofold:
\begin{itemize}
\item Introduce a general mathematical framework, and analyze the consistency, convergence and cost of the
  {proposed} algorithms in an abstract setting and in specific situations where calculations can be made
  explicit (Wishart or subgaussian distributions). The main outcome is that the subsystems drawn from the
  original system have to be chosen rectangular and not square (as initially proposed in the SADM method) to
  {yield} convergent and efficient algorithms. In other words, instead of drawing $n \times n$ subsystems, we
  will show that in many cases of applications, it is more interesting to draw $n \times n+2$ or $n \times 2n$
  subsystems in order to control the expectation {and} variance of the distribution.
\item Apply these results to revisit and improve the SADM method. This is mainly achieved in
  Section~\ref{sec:num} by considering a simple, three-point charge model of water.
\end{itemize}

\section{Mathematical setting}

{Let us} now describe more precisely the problematic.

{\subsection{General least-squares problems}}

Let $(\Omega,\mu)$ be a probability space $\Omega$ equipped with a
probability measure $\mu$. For a given arbitrary function $f \in \Lr^2(\Omega)$
and $n$ given functions $\gamma_j(x) \in \Lr^2(\Omega)$, $j =
1,\ldots,n$ all taking values in $\R$, we consider the problem of
approximating $f(x)$ by a linear combination of the functions
$\gamma_j(x)$, $j = 1,\ldots,n$.

Ideally, we would like to solve the problem of finding $\alpha = (\alpha_j)_{j = 1}^n \in \R^n$, minimizing
the function

\begin{equation}
\label{E1}
\R^n \ni a \mapsto
\Norm{f(x) -  \sum_{j = 1}^n a_j\gamma_j(x) }{\Lr^2(\Omega)}^2.
\end{equation}

The actual quality of the least-squares approximation is given by the size of
the residue $\Norm{\rho(\alpha)}{\Lr^2(\Omega)}$ where for $a =(a_j)_{j = 1}^n
\in \R^n$,

\begin{equation}
\label{Eresidu}
\rho(a)(x) = f(x) - \sum_{j = 1}^n a_j \gamma_j(x).
\end{equation}

Many minimization problems arising in mathematics and in physics can be
 {stated} under this form, for instance:
\begin{itemize}

\item[\bf{(a)}] $\Omega = [a,b]^{n}$ with two real numbers $a$ and $b
  > a$, and equipped with the measure $\dd \mu(x) = (b-a)^{-n} \dd x$
  where $\dd x$ is the Lebesgue measure on $\R^n$. Taking $\gamma:
  \Omega \to \R^{n+1}$ defined by $\gamma_i(x) = x_i$ for all
  $i\in\{1,\ldots,n\}$ and $\gamma_{n+1}\equiv 1$, the problem is
  equivalent to finding $\beta \in \R$ and $\alpha \in \R^n$
  minimizing the function
  $$
  \Norm{f(x) - \beta - \LL \alpha, x \RR}{\Lr^2([a,b]^n)}^2
  $$
  where $\LL \, \cdot\, , \cdot \, \RR$ is the standard Euclidean
  product in $\R^n$.  This is nothing  {else than a}
  multivariate linear interpolation.

  Similarly, any polynomial approximation problem in
  $\Lr^2([a,b]^n,\mu)$, where $\mu$ is a weight function,
  can be written  {in} the form \eqref{E1} by taking as $\gamma_j$ a
  basis of polynomials in dimension $n$.

\item[\bf{(b)}] Taking $\Omega = \R^n$ equipped with
  a given $n$-dimensional Gaussian measure leads to many different situations:
  The approximation by Hermite functions in $\R^n$ if $\gamma_j$ are polynomials,
  the approximation of $f$ by Gaussian chirps signal \cite{meyer-xu-97} in the
  case where $\gamma_j(x)$ are oscillating functions of $x$, or alternatively
  approximation by Gaussian wavepackets functions \cite{heller-75}
  in the context of molecular dynamics.

\item[\bf{(c)}] Consider $\Omega = \{1,\ldots,M\}$ with $M \gg n$ equipped with the
  uniform probability measure $M^{-1}\sum_{i=1}^M\delta_i$. In this
  case, an application $f$ is represented by a vector $b \in \R^M$,
   {whereas} $\gamma$ is
  represented by a matrix $A$ with $n$ columns and
  $M$ lines. The problem is then equivalent to the problem of finding
  $\alpha \in \R^n$  {that minimizes}
  $$
  \Norm{A \alpha - b }{2}^2
  $$
  where $\Norm{\cdot}{2}$ is the Euclidean norm on $\R^M$. This corresponds to the case described in
  \eqref{Eurpb}.

\item[\bf{(d)}] Consider $\Omega = \R^n \times \Omega'$ equipped with
  the measure $\mu\otimes\nu$ where $\mu$ and $\nu$ are probability
  measures on $\R^n$ and $\Omega'$ respectively.
  Taking $f(x,\omega') =  h(x) + X(\omega')$ where $X(\omega')$ is a given
  random variable on $\Omega'$, and $\gamma_j(x,\omega') = x_j$ for $j = 1,\ldots,n$ yields the
  problem of minimizing
  \begin{equation}
    \label{eq:filt}
    \min_{\alpha \in \R^n} \E\Big[\Norm{\LL \alpha , x\RR - f(x,\omega')}{\Lr^2(\R^n)}^2\Big]
  \end{equation}
  which corresponds to the linear regression of a function
  observed with some independent noise.
\end{itemize}
%

The problem \eqref{E1} is equivalent to solving the linear equation
$$
\LL \gamma, \gamma^T  \RR_{\Lr^{2}} \alpha = \LL \gamma, f \RR_{\Lr^{2}}
$$
where $\alpha = (\alpha_i)_{i = 1}^n$ and $\LL \gamma, \gamma^T \RR_{\Lr^{2}}$
is the $n \times n$ matrix with coefficients $\LL \gamma_i, \gamma_j
\RR_{\Lr^{2}}$, $i,j = 1,\ldots,n$.

If the family $(\gamma_i(x))_{i = 1}^n$ defines a full rank set of
elements of $\Lr^2(\Omega)$, the matrix $\LL \gamma, \gamma^T
\RR_{\Lr^{2}}$ is invertible, and the solution of the previous
equation reads
\begin{equation}
  \label{E2}
  \alpha = \LL \gamma, \gamma^T  \RR_{\Lr^{2}}^{-1} \cdot\LL \gamma, f \RR_{\Lr^{2}}.
\end{equation}

Apart from specific situations, where, for instance, the $\gamma_j$ can be
assumed orthogonal, the numerical approximation of \eqref{E2} is extremely
costly with respect to the dimension of $\Omega$ (see for instance \cite{bjorck-96}). Typically,
discretizations of problems of the form $\mathbf{(a)}$ yields a problem of the
form $\mathbf{(c)}$ with $m = N^n$ where $N$ is the number of interpolation
points in $[a,b]$  {needed} to approximate the
$\Lr^2$ integrals. For $n = 30$, this method is not tractable in practice, even
if $N=2$.

To avoid this {\em curse of dimensionality}, an
 {alternative would consist in approximating} the
integrals in the formula \eqref{E2} by using Monte Carlo methods. In large
dimension, the matrix $\LL \gamma, \gamma^T \RR_{\Lr^{2}} $ is,
 {however}, often ill-conditioned, and obtaining
a correct approximation of the inverse of this matrix
 {might} require in practice a very large number
of draws to minimize the error in the value of $\alpha$.

{\subsection{Principle of the algorithm}
\label{sec:algo1}}

In this abstract mathematical setting, the principle lying behind the SADM method can be extended to the
following: Retaining the idea of drawing subsystems of the original problem, we consider the following
algorithm:
\begin{itemize}
\item Draw $m$ points $X^{(i)}$, $i = 1,\ldots,m$ in $\Omega$ {independent and identically distibuted
    (i.i.d.) with distribution $\mu$.}
\item Solve the $m\times n$ least-squares sub-problem by {determining $\beta$} minimizing the function
  \begin{equation}
    \label{Esubpbs}
    \R^n \ni \beta \mapsto
    \sum_{i = 1}^m |f(X^{(i)}) -  \sum_{j = 1}^m \beta_j\gamma_j(X^{(i)})|^2.
  \end{equation}
\item Approximate the {expectation} $\bar \beta$ of the random variable $\beta$ by a Monte-Carlo method and
  analyse its distribution.
\end{itemize}


{More precisely, we define $X:=(X^{(1)},\ldots,X^{(m)})$ and the functions $F:\Omega^m \to \R^m$ and
$\Gamma:\Omega^m \to \mathcal{L}(\R^m,\R^n)$ by the {formulae}
\begin{equation}
  \label{E3}
  \forall i = 1,\ldots,m, \quad
  F_i(x^{(1)},\ldots,x^{(m)}) = f(x^{(i)})
\end{equation}
and
\begin{equation}
  \label{E4}
  \forall i = 1,\ldots,m,  \quad\forall\, j =1 ,\ldots,n, \quad
  \Gamma_{ij}(x^{(1)},\ldots,x^{(m)}) = \gamma_j(x^{(i)}).
\end{equation}}
{The random vector $\beta$ then minimizes the function
$$
\beta \mapsto \Norm{F(X) - \Gamma(X) \beta}{2}^2,
$$
where $\Norm{\cdot}{2}$ is the standard Euclidean norm on $\R^m$.}

{Under the assumption that $\Gamma^T(X)\Gamma(X)$ is invertible almost surely (a.s.),
\begin{equation}
  \label{E5}
  \beta = R(X) F(X) := ((\Gamma^T \Gamma)^{-1} \Gamma^T)(X) F(X).
\end{equation}
The expectation of $\beta$ is then given by the formula
\begin{equation}
  \label{E6}
  \bar{\beta}:=\E \beta = \int_{\Omega^m} ( (\Gamma^T \Gamma)^{-1} \Gamma^T
  F)(x^{(1)},\ldots,x^{(m)})\,
  \dd \mu(x^{(1)})\otimes \cdots\otimes \dd\mu(x^{(m)}).
\end{equation}
Our algorithm consists in using the Monte-Carlo method to compute the previous expectation: we approximate
$\bar\beta$ by
\begin{equation}
  \label{eq:est}
  \bar{\beta}_N=\frac{1}{N}\sum_{i=1}^N\beta_{i},
\end{equation}
where $\beta_{i}, i\geq 1$ are i.i.d.\ realizations of the random vector
$\beta\in\R^n$, obtained by~(\ref{E5}) from i.i.d.\ realizations of
the random $n\times m$ matrix $\Gamma(X)$.}

\bigskip

{Of course, one expects that $\bar{\beta}$ should converge to the solution of the least square
problem \eqref{E2} when $m\rightarrow+\infty$. This indeed can be easily proved under the additional
assumption that $f$ and $\gamma_j$, $1\leq j\leq n$ belong to $\Lr^2(\Omega)$. By the strong law of large
numbers,
\begin{equation}
  \label{eq:MC1}
  \frac{1}{m}(\Gamma(X)^T\Gamma(X))_{ij}
  =\frac{1}{m}\sum_{k=1}^m\gamma_i(X^{(k)})\gamma_j(X^{(k)})
\end{equation}
converges $\P$-a.s.\ to $(\LL \gamma, \gamma^T \RR_{\Lr^{2}})_{ij}$ when
$m\rightarrow+\infty$. Similarly,
\begin{equation}
  \label{eq:MC2}
  \frac{1}{m}(\Gamma(X)^TF(X))_{i}
  =\frac{1}{m}\sum_{k=1}^m\gamma_i(X^{(k)})f(X^{(k)})
\end{equation}
converges $\P$-a.s.\ to $(\LL \gamma, f \RR_{\Lr^{2}})_i$. Consequently, if the matrix $\LL \gamma,
\gamma^T \RR_{\Lr^{2}}$ is invertible,
\begin{equation}
  \label{eq:expr}
  \beta=\Big(\frac{1}{m}\Gamma^T(X)\Gamma(X)\Big)^{-1}\frac{1}{m}\Gamma(X)^TF(X)
\end{equation}
converges $\P$-a.s.\ to $\alpha$ given by \eqref{E2} when $m\rightarrow+\infty$.}
\medskip

{However our goal is not to analyse more finely this convergence, as we are concerned with situations where
  the least square problem~\eqref{E1} is ill-posed or computationally unfeasible due to the high diemsnion of
  the problem. In the opposite, considering the case where $m$ is {\em small} in comparison with the dimension
  of $\Omega$ ($M$ in the case of SADM) should reduce the computational cost, provided that the efficiency of
  the Monte-Carlo approximation is good.} To express the fact that we are in a regime where $m$ is small, {\bf
  we assume in the following that $\boldsymbol{m \leq Cn}$ }for some constant $C$ (typically $m = n+2$ or $m =
2n$ for practical applications).

Therefore, to make sure that the previous algorithm is efficient, we have to verify the following points: 

\begin{itemize}
\item[(i)] The random variable $\beta$ has {finite expectation} and variance. Here the bounds may depend on
  $n$, but not on the cardinal of $\Omega$ ($M$ in the SADM description above). This condition is crucial to
  ensure the convergence of a Monte-Carlo method and the approximability of $\bar \beta$. {In addition, the
    smaller is the variance, the faster the Monte-Carlo approximation converges to $\bar{\beta}$.}
\item[(ii)] The average $\bar \beta$ is a good alternative to the solution of the original problem \eqref{E1}
  in the sense that $\bar \beta - \alpha = \mathcal{O}(\Norm{\rho(a)}{})$ where $\rho(a)$ is the residue
  \eqref{Eresidu}. In other words, if $f$ is close to a linear combination of the functions $\gamma$ the
  residue will be small and the standard least-square approximation will be efficient. In this situation,
  $\bar \beta $ will also lead to a good approximation, and be close to the solution $\alpha$. On the other
  hand, when the residue is large, $\bar \beta$ and $\alpha$ may differ, but in this situation the
  approximation of $f$ by a linear combination of the functions of $\gamma$ is poor in any case.
\end{itemize}

{In Section~\ref{sec:sigma0} we give various conditions that warrant the latter requirements. In
particular, we study the consistency of the algorithm, give conditions ensuring the convergence of the Monte
Carlo method, and analyze the computational cost. In the specific instance where $\Gamma(X)^T\Gamma(X)$ has
the Wishart distribution, all computations can be made explicitly, and we obtain precise estimates and an
optimal choice of the parameter $m$. The two values $m=n+2$ and $m=2n$ are of specific interest in this
situation. In addition, we prove that the choice $m=n$} leads to a random variable $\beta$ with infinite
expectation, which partly explains the Cauchy-like distributions observed in \cite{chipot-98} {with} the SADM
method.

{\subsection{The algorithm in the non-invertible case}
\label{sec:algo2}}

In practice, the {almost sure} invertibility of $\Gamma^T(X)\Gamma(X)$ cannot be guaranteed --- and
obviously not for problems of the form $\mathbf{(c)}$, where {all the random variables
$X^{(1)},\ldots,X^{(m)}$ may be equal with positive probability.}

In a more general setting, we,  {hence, restrict
ourselves} to realizations of $X$, such that matrix $\Gamma(X)$ is sufficiently
 {well conditioned}, in the following sense: Denoting by 
$s_1(\Gamma(X))$  the smallest eigenvalue of the symmetric positive
matrix $\Gamma(X)^T\Gamma(X)$, we only consider realizations of $X$,
such that $s_1(\Gamma(X))$ is greater than some threshold $\sigma$, which may
depend on $n$ and $m$. In this case, rather than approximating \eqref{E6}, we
will estimate the conditional expectation
\begin{equation}
\label{Esbeta}
\bar\beta^\sigma :=
\E^\sigma \beta = \E[\,\beta\mid s_1(\Gamma(X))>\sigma\,]
\end{equation}
by
\begin{equation}
  \label{eq:est-sigma}
  \bar{\beta}_N^\sigma=\frac{1}{N}\sum_{i=1}^N\beta_i^\sigma,
\end{equation}
where the $\beta^\sigma_i$ are obtained from a sequence of i.i.d.\
realizations of the random vector $\beta\in\R^n$ in~(\ref{E5}), from
which have been removed all realizations such that $s_1(\Gamma(X))\leq
\sigma$. Note that~(\ref{eq:est}) is a particular case
of~(\ref{eq:est-sigma}) for $\sigma=0$, provided that
$\P(s_1(\Gamma(X))=0)=0$.
\medskip

{Again,} such a method will be of interest in terms of computational cost if $m$ is
 {on} the order of magnitude of $n$ (in all the
applications  {considered herein}, $m = n+2$ or
$m = 2n$ will be sufficient) and if $\P(s_1(\Gamma(X))>\sigma)$ is not too
small --- because drawing a realization of $X$ such that
$s_1(\Gamma(X))>\sigma$ requires {an} average number
$\P(s_1(\Gamma(X))>\sigma)^{-1}$ of realizations of $X$.

From  {the perspective of precision}, this method
will  {perform well} if the variance of $\beta$
conditionally on $\{s_1(\Gamma(X))>\sigma\}$ has an appropriate behavior with
respect to $n$ and $m$, and if $\bar\beta^\sigma$ defined in \eqref{Esbeta}
provides a good approximation of the solution of the original least-squares
problem.

{The specific case where the $\Gamma(X)^T\Gamma(X)$ is not a.s.\ invertible is studied in
Section~\ref{sec:sigma}, where} we give various conditions that warrant the latter requirements. {The
instance where $\Gamma(X)$ has subgaussian entries (which covers the Wishart case mentionned above) is then
studied in more details and leads again to optimal choices of $m$, $N$ and $\sigma$.}


\bigskip

\section{The invertible case}
\label{sec:sigma0}

{In all this section, we assume that the matrix $\Gamma(X)^T\Gamma(X)$ is a.s.\ invertible.

\subsection{Preliminary results}

Before studying the algorithm of Section~\ref{sec:algo1}, let us define for $q\in[2,+\infty]$
\begin{equation}
\label{EKmq}
K_q(\Gamma):= \left[ \E\,  \frac{1}{s_1(\Gamma(X))^{\frac{q}2}}\right]^{\frac{2}q} ,
\end{equation}
where $\Gamma(X)$ is the random matrix defined by \eqref{E4} and with the usual convention that
$K_\infty(\Gamma)=\|s_1(\Gamma(X))^{-1}\|_{\Lr^\infty}$.}  Note that $K_q(\Gamma)$ depends on $n$ and $m$.

{The proof of the next lemma is given in Appendix~\ref{sec:lemma}.}
\begin{lemma}
\label{Lborne}
Let {$p \in [1,\infty]$} and $g \in \Lr^p(\Omega)$. Let us define the function $G$ from $g$ as $F$
is defined from $f$ in \eqref{E3}. Let also $R(X)$ be the random matrix defined in \eqref{E5}.
\begin{description}
\item[\textmd{(a)}] Assume that $K_q(\Gamma) < +\infty$ where $q$ is such that $q^{-1} + p^{-1} = 1$.  Then we
  have
  \begin{equation}
    \label{Esecondest}
    \E \Norm{R(X) G(X) }{2} \leq \sqrt{n}m \sqrt{K_{q}(\Gamma)}
    \Norm{g}{\Lr^p(\Omega)}.
  \end{equation}
\item[\textmd{(b)}] Assume that {$p \in [2,\infty]$} and that
  $K_q(\Gamma) < +\infty$ where $q$ is such that $2q^{-1} + 2p^{-1} =
  1$.  Then we have
  \begin{equation}
    \label{Efirstest}
    \E \Norm{R(X) G(X) }{2}^2 \leq nm^{2} K_{q}(\Gamma) \Norm{g}{\Lr^p(\Omega)}^2.
  \end{equation}
\end{description}
\end{lemma}

{The next result is a first consequence of this lemma. {We recall} the definition of $\bar\beta$ in
\eqref{E6} and that $\rho(a)$ denotes the residue \eqref{Eresidu} associated with the function $f$ and the
coefficients $a_j$, $j = 1,\ldots,n$.
\begin{proposition}
  \label{prop:consist}
  Let $a = (a_j)_{j = 1}^n \in \R^n$ and $m \leq Cn$ for some constant $C$.  Assume that $\rho(a) \in
  \Lr^p(\Omega)$ {and $K_q(\Gamma) < +\infty$ for some $p\in[1,+\infty]$ and with $q^{-1}+p^{-1}=1$.}  Then
  there exists a constant $C(n)$ {depending on $n$} such that
  \begin{equation}
    \label{Econsist}
    \E\Norm{\beta-a}{2}\leq C(n)  \Norm{\rho(a)}{\Lr^p(\Omega)}.
  \end{equation}
\end{proposition}
}

\begin{Proof}
{By definition of $R(X)$, and as $\Gamma(X)^T \Gamma(X)$ is invertible, we have 
$$
R(X) \Gamma(X) a = a. 
$$
Hence 
\begin{equation}
\label{Econss}
\beta - a = R(X)F(X) - R(X)\Gamma(X) a = R(X) \rho(a)(X).
\end{equation}
where $\rho(a)(X)$ is defined from $\rho(a)$ as $F$ was defined from $f$ in \eqref{E3}. }
{The result then follows from Lemma~\ref{Lborne}~(a) with $C(n) = \sqrt{n}m \sqrt{K_q(\Gamma)}$.}
\end{Proof}

\subsection{Average and variance}

The following result {is an immediate consequence of Prop.~\ref{prop:consist}. It} gives conditions on $f$ and
$\Gamma$ to ensure that the random variable $\beta$ has finite expectation, {and thus that the Monte Carlo
approximation a.s.\ converges to $\bar{\beta}$ when $N\rightarrow+\infty$.}
\begin{corollary}
\label{cor:average}
Let $m \leq Cn$ for some constant $C$ and assume that $f \in \Lr^p(\Omega)$ {and $K_q(\Gamma) < +\infty$
for some $p \in [1,+\infty]$ and with $q^{-1}+p^{-1}=1$.}  Then there exist a constant $C(n)$ depending on
$n$ such that
$$
\E \Norm{\beta}{2} \leq C(n) \Norm{f}{\Lr^p(\Omega)}. 
$$
\end{corollary}



\bigskip

In order to estimate the convergence rate of algorithm, we need to construct confidence regions with
asymptotic level (less than) $\eta$ for the Monte Carlo approximation of $\bar{\beta}$. We are going to
consider confidence regions of the form $[a_1,b_1]\times\ldots\times[a_n,b_n]$, by taking each $[a_i,b_i]$ as
a confidence interval of asymptotic level $\eta/n$ for the $i$-th coordinate $\beta_i$ of $\beta$. Note that
more precise asymptotic confidence regions exist --- see {for instance}~\cite{anderson-84} --- but the
previous confidence region is more convenient for computation. Note also that non-asymptotic estimates could
be obtained using Berry-Essen-type inequalities --- see for instance \cite{petrov-75}.

{This leads to the choice}
$$
b_i-a_i=2x(n,\eta)\sqrt{\mbox{Var}(\beta_i)/N},\quad\forall i\in\{1,\ldots,n\}
$$
{where $N$ is the number of draws in \eqref{eq:est}, and}
where $x(n,\eta)>0$ is {the solution of}
\begin{equation}
  \label{eq:def-x}
  \frac{1}{\sqrt{2\pi}}\int_{x(n,\eta)}^{+\infty}e^{-u^2/2}du=\frac{\eta}{2n}.
\end{equation}
In this case, the {Euclidean} diameter of the confidence region is bounded by
\begin{equation}
  \label{eq:conf2}
  2x(n,\eta)\sqrt{\mbox{Tr}(\mbox{Cov}(\beta))/N},
\end{equation}
{where $\text{Cov}(\beta)$ is the covariance matrix of the random vector $\beta$, defined by
$$
\mbox{Cov}(\beta):=\E[(\beta-\E \beta)(\beta-\E \beta)^T].
$$

The next result gives bounds on the quantity $\mbox{Tr}(\mbox{Cov}(\beta))$, which, in view
of~\eqref{eq:conf2}, controls the rate of convergence of the Monte-Carlo approximation.}

\begin{proposition}
  \label{Tmoments}
  Let $m \leq Cn$ for some constant $C$ and assume that $\rho(\bar{\beta})\in\Lr^p(\Omega)$ {and
    $K_q(\Gamma) < +\infty$ for some $p\in[2,+\infty]$ and with} $2p^{-1}+2q^{-1}=1$. Then there exist a
  constant $C(n)$ depending on $n$, such that
   \begin{equation}
    \label{eq:var*1}
    \mbox{\textup{Tr}}(\mbox{\textup{Cov}}(\beta))\leq C(n) 
    \Norm{\rho(\bar\beta)}{\Lr^p(\Omega)}^2. 
  \end{equation}
\end{proposition}

\begin{Proof}
  {Let} $g = \rho(\bar\beta)$ and {define $G$ from $g$ as $F$ is defined from $f$} by \eqref{E3}. {Then}
  $$
  \mbox{Tr}(\mbox{Cov}(\beta))
  =\E\Norm{\beta-\bar{\beta}}{2}^2
  = \E\Norm{R(X) G(X)}{2}^2.
  $$
  The result, {hence}, follows from Lemma~\ref{Lborne}~(b) with $C(n) = n m^2 K_{q}(\Gamma)$. 
\end{Proof}


{These results show that the convergence {of our algorithm} relies on an assumption of the form
  $K_q(\Gamma)<+\infty$, which corresponds to the finiteness of a negative moment of the random variable
  $s_1(\Gamma(X))$. Such an assumption is clearly problem-dependent and has to be checked in each specific
  problem considered. Conditions ensuring this property when $q<+\infty$ are given in
  Appendix~\ref{sec:moment} and are used to handle the specific case of Wishart matrices in
  Subsection~\ref{sec:Wishart}.}

{Note that, under the assumptions of this section, the condition $K_q(\Gamma)$ is unlikely to be satisfied
  when $q=+\infty$. Indeed, since $\Gamma^T(X)\Gamma(X)$ is assumed a.s.\ invertible, the measure $\mu$ must
  have no atom, and hence $\Omega$ is continuous (\emph{i.e.} not denumerable). If we assume in addition that
  the functions $\gamma_j$ are regular on $\Omega$, so are the eigenvalues of $\Gamma^T(x)\Gamma(x)$ as a
  function of $x=(x^{(1)},\ldots,x^{(m)})\in\Omega^m$. Since the smallest eigenvalue is $0$ when
  $x^{(1)}=\ldots=x^{(m)}$, we deduce that $\mathbb{P}(s_1(\Gamma(X))<\eta)>0$ for all $\eta>0$, which means
  that $K_\infty(\Gamma)=\infty$. The way to handle the case $q=\infty$ is explained in
  Section~\ref{sec:sigma}.}

\subsection{Link with the least square approximation}

Formula~\eqref{E6} proposes an alternative solution $\bar{\beta}$ to the solution $\alpha$ given by~(\ref{E2})
of the least-squares problem \eqref{E1}. We now {provide} estimates between these two solutions.

A precise error estimate depends on the tackled problem {(see for instance Section~\ref{sec:Wishart}).}
{Here}, we give a general result.

{\begin{proposition}
  \label{prop:comp-lsquares}
  Assume that $f,\gamma_1,\ldots,\gamma_n$ belong to $\Lr^2(\Omega)$ and that $K_2(\Gamma)<+\infty$. Then
  there exists a constant $C(n)$ such that 
  \begin{equation}
    \label{Econsist2}
    \Norm{\bar\beta - \alpha}{2} \leq C(n)  \Norm{\rho(\alpha)}{\Lr^2(\Omega)}.
  \end{equation}
\end{proposition}

\begin{Proof}
Observing that $\rho(\alpha)\in\Lr^2(\Omega)$, this is an immediate consequence of Prop.~\ref{prop:consist}
and of the inequality $\Norm{\bar{\beta}-\alpha}{2}\leq\E\Norm{\beta-\alpha}{2}$.
\end{Proof}}

{In other words, the better $f$ can be approximated by a linear combination of the functions $\gamma_j$,
$1\leq j\leq n$, the closer the result of our algorithm is from the actual least square approximation.}


\subsection{ {Computational} cost of the {algorithm}}
\label{sec:cost}

Let $\varepsilon$ be a required precision for the approximation of
$\bar\beta = \E \beta$ by the Monte Carlo simulation
\eqref{eq:est}. For large $N$, using \eqref{eq:conf2}, we must take
$$
N \sim 4 x(n,\eta)^2 \varepsilon^{-2}  \mbox{Tr}(\mbox{Cov}(\beta)).
$$
Since, for all $x>0$,
\begin{equation}
  \label{eq:ineq-gauss}
  \int_{x}^{+\infty}e^{-u^2/2}du\leq\frac{1}{x}\int_{x}^{+\infty}ue^{-u^2/2}du
  =\frac{e^{-x^2/2}}{x},
\end{equation}
{we deduce from~\eqref{eq:def-x} that, for $n/\eta$ large enough,}
\begin{equation}
  \label{eq:conf2-cst}
  x(n,\eta)^2\leq\Big(2\log\frac{n\sqrt{2}}{\eta\sqrt{\pi}}\Big)
\end{equation}

In addition, each step of the algorithm requires {to evaluate the} matrix $\Gamma(X)^T\Gamma(X)$ and the
vector $\Gamma(X)^T F(X)$ and {to invert} the matrix $\Gamma(X)^T\Gamma(X)$. The cost of
these {operations} is of order $C n^2 m $.

{Hence,} we see that the cost of the algorithm is of order
$$
C \varepsilon^{-2} m n^2 \log n \,\mbox{Tr}(\mbox{Cov}(\beta)).
$$
Under the hypothesis of Proposition \ref{Tmoments} {and using the explicit expression of $C(n)$ obtained in
  the proof of this proposition,} the computational cost {can be bounded by}
\begin{equation}
\label{EKost}
C \varepsilon^{-2} m^3 n^3 \log n \,  K_{q}(\Gamma) \Norm{\rho(\bar\beta)}{\Lr^p(\Omega)}^2
\end{equation}
for $2p^{-1}+2q^{-1}=1$.

{It may be observed }that this cost depends only on $n$ and $m$ --- and not the dimension of
$\Omega$. Moreover, it depends on the least-squares residue of the problem \eqref{E1}. In the {event} where
$f$ is close to a linear combination of the functions $\gamma_j$, the algorithm is, {therefore}, cheaper
{(and, by Prop.~\ref{prop:comp-lsquares}, more precise).}  {As a consequence}, the cost of our algorithm is
driven by the quality of the original least-squares approximation in Problem \eqref{E1}.


\subsection{The Wishart case}
\label{sec:Wishart}

{Let us now} consider the case where
$\Omega=\R^n$,
$$
\dd\mu(x)=(2\pi)^{-n/2}\exp(-\Norm{x}{2}^2/2)\dd x_1\ldots \dd x_n
$$
and $\gamma_j(x)=x_j$ for $j\in\{1,\ldots,n\}$ ---
 {{\it i.e.}} linear interpolation. In this case, {the}
random vectors $X^{(i)}$ are standard $n$-dimensional Gaussian vectors, the
matrix $\Gamma(X)$ is a $m \times n$ matrix with i.i.d.\ standard Gaussian
entries and the law of the matrix $\Gamma(X)^T\Gamma(X)$ is the so-called
Wishart distribution --- see {\it e.g.}~\cite{anderson-84}.

The joint distribution of its eigenvalues is known explicitly and can be found for example
in~\cite[p.534]{anderson-84}. In particular, $\Gamma(X)^T\Gamma(X)$ is a.s.\ invertible {if} $m\geq n$. The
explicit density of the eigenvalues has been used to {obtain} estimates on the law of the smallest eigenvalue
of such matrices in~\cite{edelman-88,edelman-91,chen-dongarra-05}. {These results allow us to obtain
  explicit estimates in the Wishart case, proved in Appendix~\ref{sec:Wishart-pf}. We shall restrict here to
  the case where $f$ and $\rho(\bar{\beta})$ belong to $\Lr^\infty(\Omega)$, and we refer to
  Appendix~\ref{sec:Wishart-pf} for further estimates.}

{Under the previous assumptions, the conditions of Corollary~\ref{cor:average} and
  Proposition~\ref{Tmoments} are satisfied for all $m\geq n+2$. The computational cost is (asymptotically)
  minimal for the choice $m=2n$ and the corresponding computational cost is bounded by
\begin{equation}
  \label{eq:cost-Wishart}
  C\varepsilon^{-2} n^5\log n  \Norm{\rho(\bar\beta)}{\Lr^\infty}^2  
\end{equation}
for an explicit constant $C$ independent of $n$, where $\varepsilon$ is the required precision of the
algorithm.

In addition, the consistency error of Proposition~\ref{prop:comp-lsquares} is bounded by
$$
C'n\Norm{\rho(\alpha)}{\Lr^\infty}
$$
for a constant $C'$ independent of $n$.

We, hence, see that the values $m=n+2$ and $m=2n$ are of specific interest in terms of convergence and
computational cost. Although the Wishart case corresponds to very simple approximation problems, this result
gives valuable clues about the way parameters should be chosen in our algorithm. These specific values of $m$
are numerically tested in the example of the three-point charge model of water developed in
Section~\ref{sec:num}, where improvements of the SADM method are considered.}

\section{The general case}
\label{sec:sigma}

 {Let us now consider} the general case where
$\Gamma(X)^T \Gamma(X)$ is not assumed to be a.s. invertible.

{Fix $\sigma > 0$.} We denote by $\E^\sigma$ (resp.
$\mbox{\textup{Cov}}^{\sigma}$) the expectation (resp.\ covariance matrix)
conditionally on the event $\{s_1(\Gamma(X))>\sigma\}$. As an approximation of
the solution of the  {least-squares} problem, we
will  {examine} the conditional expectation
\begin{equation}
\label{Econdbeta}
\bar\beta^{\sigma} = \E^{\sigma}(\beta).
\end{equation}

{As will appear below, our algorithm always converges for any $\sigma>0$. As in the invertible case, its
  performance relies on precise estimates on convergence, consistency and computational cost, given
  below. Optimal computations will then be detailed in the specific instance where the matrix $\Gamma(X)$
  has independent sub-Gaussian entries.}

\subsection{Consistency, convergence and  {computational} cost}

{We first generalize Proposition~\ref{prop:consist}: For all $q\in[1,+\infty]$, let
\begin{equation}
  \label{eq:def-K-sigma}
  K^\sigma_q(\Gamma):=\left[\E^\sigma\frac{1}{s_1(\Gamma(X))^{\frac{q}{2}}}\right]^{\frac{2}{q}}.    
\end{equation}}
\begin{proposition}
  \label{prop:consist2}
  Let $a_j$, $j = 1,\ldots,n$ be $n$ numbers $a_j$. Assume that $\rho(a) \in\Lr^p(\Omega)$ for some
  {$p\in[1,+\infty]$, then
  $$
  \Norm{\bar\beta^{\sigma}-a}{2}\leq \E^\sigma\Norm{\beta-a}{2}\leq
  \frac{\sqrt{n}m}{\P(s_1(\Gamma(X))\geq\sigma)^{1/p}}\,\sqrt{K^\sigma_q(\Gamma)}\,
  \Norm{\rho(a)}{\Lr^p(\Omega)}
  $$}
  where $q$ is such that $q^{-1} + p^{-1} = 1$.
\end{proposition}

\begin{Proof}
  {Using the inequality}
  $$
  \E^\sigma\Norm{\rho(a)(X)}{1}^p
  \leq\frac{\E\Norm{\rho(a)(X)}{1}^p}{\P(s_1(\Gamma(X))\geq\sigma)}
  $$
  {in~\eqref{eq:pf-lemma},} the proof is exactly the same as {that put forth} in Lemma \ref{Lborne} and
  Proposition~\ref{prop:consist}.
\end{Proof}

{Note that, by definition of $\E^\sigma$, for all $q\in[1,\infty]$,
\begin{equation}
  \label{eq:trivial-ineq}
  K^\sigma_q(\Gamma)\leq\sigma^{-1}.    
\end{equation}
}
{In particular, taking $a=0$ in} {the} previous result implies that conditional
expectation~(\ref{Econdbeta}) is always well defined for $\sigma>0$ as soon as $f\in\Lr^p(\Omega)$ for some
$p\in[1,+\infty]$.

The following result generalizes Proposition \ref{Tmoments} to the case where $\sigma > 0$. Its proof is very
similar to {that} of Proposition \ref{Tmoments}. We {will, hence,} omit it {here}.
\begin{proposition}
  \label{Tmoments2}

  Assume that the function $\rho(\bar\beta^{\sigma}) \in \Lr^p(\Omega)$ for
  {$p\in[2,+\infty]$. We have
  \begin{equation}
    \label{eq:var1}
    \mbox{\textup{Tr}}(\mbox{\textup{Cov}}^{\sigma}(\beta))\leq
    \frac{nm^2}{\P(s_1(\Gamma(X))\geq\sigma)^{2/p}}\,K^\sigma_q(\Gamma)\,
    \Norm{\rho(\bar\beta^{\sigma}) }{\Lr^p(\Omega)}^2
  \end{equation}
  where $q$ is such that $2q^{-1}+2p^{-1}=1$.}
\end{proposition}

{Although the trivial inequality~\eqref{eq:trivial-ineq} always allows one to {infer} explicit bounds from
  the previous results, there are cases where optimal estimates on $K^\sigma_q(\Gamma)$ are much better. Since
  our performance analysis relies heavily on precise estimates on $K^\sigma_q(\Gamma)$, it is desirable to
  obtain conditions for better estimates. Such conditions are given in Proposition~\ref{prop:mom2} in
  Appendix~\ref{sec:moment}, and will be used to handle the sub-Gaussian case described in the next
  subsection.}

\bigskip

We now consider the cost of the algorithm: Let $\beta^{(\sigma)}$
denote a random variable having the law of $\beta$ conditioned on
$\{s_1(\Gamma(X))\geq\sigma\}$. The cost of the algorithm is
determined by
\begin{itemize}
\item the number $N$ of simulations of $\beta^{(\sigma)}$ needed to
  ensure that the diameter of the confidence region for the
  Monte Carlo estimation of
  $\E(\beta^{(\sigma)})=\E^{\sigma}(\beta) = \bar\beta^\sigma$ is
  smaller than a given precision $\varepsilon$. To control this, we
  use the upper bound on the confidence region diameter given
  by~(\ref{eq:conf2}), where $\eta$ is the level of confidence of the
  approximation;
\item the average number of draws of the random variable $X$ needed to
  simulate a realization of $\beta^{(\sigma)}$, which is
  {$1/\P(s_1(\Gamma(X))\geq\sigma)$.} Note that a draw corresponds
  to simulating a $nm$-dimensional random variable.
\item the computation of the $n\times n$ matrix
  $\Gamma(X)\Gamma(X)^T$, which is of order $n^2 m$ --- all other
  computational costs, including the cost of the computation of
  $s_1(\Gamma(X))$ or the inversion of $\Gamma(X)^T\Gamma(X)$, are of
  a smaller order with respect to the dimension $n$ of the problem,
  provided  {that} $m \geq n$.
\end{itemize}

Consequently, the cost of the algorithm is bounded by
$$
 C N\,
\P(s_1(\Gamma(X))\geq\sigma)^{-1}(nm + n^2m)
$$
for some constant $C>0$. As
$$
N \sim 4 x(n,\eta)^2 \varepsilon^{-2}\mbox{Tr}(\mbox{Cov}^{\sigma}(\beta)),
$$
because of~(\ref{eq:conf2-cst}), the cost can be bounded by
$$
C  \varepsilon^{-2}n^2m\log n\,
\P(s_1(\Gamma(X))\geq\sigma)^{-1}
\mbox{Tr}(\mbox{Cov}^{\sigma}(\beta)).
$$
{Thus}, if
$\rho(\bar{\beta}^\sigma)\in\Lr^p(\Omega)$ for {$p\in[2,+\infty]$,} because of
Proposition~\ref{Tmoments2}, the cost is bounded by
{$$
C\varepsilon^{-2}n^3m^3\log n\, \P(s_n(\Gamma(X))\geq\sigma)^{-1-\frac{2}{p}}\,K_q^\sigma(\Gamma)\,
\Norm{\rho(\bar\beta^{\sigma}) }{\Lr^p(\Omega)}^2
$$}
for some constant $C>0$, {where $2q^{-1}+2p^{-1}=1$.}

 {We, hence, can} see that the choice of an
optimal threshold $\sigma$ has to be balanced to optimize the ratio between
{$K^\sigma_q(\Gamma)$ and the} probability $\P(s_n(\Gamma(X))\geq\sigma)$
at some appropriate powers.
\medskip

{Again, explicit bounds may depend on the tackled problem. Hereafter, we develop the particular instance}
where $\Gamma(X)$ is a matrix with independent sub-Gaussian entries.

\subsection{The sub-Gaussian case}
\label{sec:subgaussian}

{We recall that the convergence of the algorithm holds for any choice of $\sigma>0$. The goal of this section
  is to study the behaviour of the computational cost in the subgaussian case as a function of $\sigma$ and
  $m$.}

{We consider} the case where $\Omega=\R^n$,
$$
\dd\mu(x)=\otimes_{i=1}^n\dd\nu(x_i)
$$
for some probability measure $\nu$ on $\R$, and $\gamma_j(x)=h(x_j)$ for $j\in\{1,\ldots,n\}$ for some
function $h$ on $\R$. This {is tantamount} to the case of an approximation of the function $f$ on $\R^n$ by a
linear combination of functions depending {on only} one variable.

In this case, it is clear that all the entries of matrix $\Gamma(X)$ are i.i.d.  Let us assume that these
random variables are \emph{sub-Gaussian}, \emph{i.e.}
$$
\forall t>0,\quad \nu(\{x\in\R:|h(x)|>t\})\leq 2\exp(-t^2/R^2)
$$
for some $R>0$. {Such} is the case, in particular if $h$ is bounded or if $\nu$ has compact support and $h$
is continuous on the support of $\nu$. {{\sc Rudelson} \& {\sc Vershynin}~\cite{rudelson-vershynin-08} have
recently obtained estimates on the distribution of $s_1(\Gamma(X))$ in the subgaussian case, optimal in the
sense that they are consistent with the explicit bounds in the Wishart case. 

{Using these results, under the assumption that $f$ and $\rho(\bar{\beta})$ belong to
$\Lr^\infty(\Omega)$ and taking $\sigma=a n$ for some constant $a>0$, computations in
Appendix~\ref{sec:sub-gaussian-pf} prove that the optimal choice for $m$ in terms of (asymptotic)
computational cost is $m=2n$, and we have the same estimates on the computational cost and the consistency as
in Section~\ref{sec:Wishart}.}

{This shows that, choosing conveniently $\sigma$, the computational cost has the same behaviour as is the
  Wishart case. In addition,} the result in terms of computational cost in $n$ appears to be relatively
unaffected by the choice of $\sigma$.
{In particular, the specific value of the constant $a$ such that $\sigma=a n$ only has an influence of the
  constant $C$ in~\eqref{eq:cost-Wishart}.}

%
\section{Improvement of the SADM method}

\label{sec:num}

The statistical analysis of distributed multipoles (SADM) algorithm
 {put forth} in \cite{chipot-98} corresponds to a
problem of the form \textbf{(c)}, where $(\alpha_j)_{j = 1}^n$ represent the
unknown multipoles borne by the $n$ particles, and $\gamma_j(x)= 1/\| x - x_j\|$ the
electrostatic potential functions, where
$x_1,\ldots,x_n$ denote the positions of the particles. The space $\Omega$ is
made of $M$ points in the three-dimensional
 {Cartesian space, lying away from the atomic
positions,} with $M >> n$.

 {However more computationally intensive than the
least-squares scheme, this pictorial approach provides a valuable information
as to whether the atomic multipoles are appropriately defined, depending on how
spread out the corresponding distributions are.} For instance, description of
the molecular electrostatic potential of dichlorodifluoromethane (CCl$_2$F$_2$)
by means of a simple point-charge model yields a counterintuitive
C$^{\delta-}$---X$^{\delta+}$ bond polarity --- where X = Cl or F, blatantly
violating the accepted rules of electronegativity differences. Whereas the
least-squares route merely supplies crude values of the charge borne by the
participating atoms, the SADM method offers a diagnosis of pathological
scenarios, like that of dichlorodifluoromethane. In the latter example, the
charge centered on the carbon atom is indeterminate, as mirrored by its
markedly spread distribution~\cite{chipot-98}. The crucial issue of buried
atoms illustrated here in the particular instance of CCl$_2$F$_2$ can be
tackled by enforcing artificially the correct bond polarity by means of
hyperbolic restraints~\cite{Bayly1993}. Violations of the classical rules of
electronegativity differences may, however, often reflect the incompleteness of
the electrostatic model --- {\em e.g.} describing an atomic quadrupole by a
mere point charge. Addition of atomic dipoles to the rudimentary point-charge
model restores the expected, intuitive C$^{\delta+}$---X$^{\delta-}$ bond
polarity~\cite{chipot-98}.

In this section, we revisit the prototypical example of the three-point charge
model of water. The molecular geometry was optimized at the MP2/6-311++G($d,p$)
level of approximation. The electrostatic potential was subsequently mapped on
a grid of 2,106 points surrounding the molecule, at the same level of theory,
including inner-shell orbitals. All the calculations were carried out with the
{\sc Gaussian 03} suite of programs~\cite{Frisch2003}. Brute-force solution of
the least-squares problem \eqref{E1}, employing the {\sc Opep}
code~\cite{Angyan2003}, yields a net charge of $-$0.782 electron-charge unit
(e.c.u.) on the oxygen atom --- hence, a charge of $+$0.391 e.c.u. borne by the
two hydrogen atoms, with a root-mean square deviation between the point-charge
model regenerated and the quantum-mechanical electrostatic potential of 1.09
atomic units, and a mean signed error of 51.1~\%. This notoriously large error
reflects the incompleteness of the model
--- a simple point charge assigned to the oxygen atom being obviously unable
to describe in a satisfactory fashion the large quadrupole borne by the latter.

On account of the ${\cal C}_{2v}$ space-group symmetry of water, only one net
atomic charge would, in principle, need to be determined --- the point charges
borne by the two hydrogen atoms being inferred from that of the oxygen atom.
Inasmuch as the SADM scheme is concerned, this symmetry relationship translates
to a single equation to be solved per realization or experiment. Without loss
of generality, two independent parameters will, however, be derived from the
electrostatic potential, the point charges borne by the two hydrogen being
assumed to be equal. Furthermore, in lieu of solving the individual
C$_{2,106}^2$ systems of 2 $\times$ 2 linear equations, incommensurable with
the available computational resources, it was chosen to select randomly 500,000
such systems.

The running averages of the charge borne by the oxygen atom are shown in
Figure~{\ref{averages}} as a function of the number of individual realizations,
for the SADM algorithm with $n = N_s$ points and its proposed enhancement,
using 2, 4 and 8 additional grid points per realization --- with the notations
 {utilized in} the previous section,
 {the latter translates to} $m = N_s +2$, $N_s
+4$ and $N_s +8$. From the onset, it can be seen that the SADM scheme yields
the worst agreement with the target value derived from the least-squares
problem \eqref{E1}, and that inclusion of supplementary equations to the SADM
algorithm rapidly improves the accord. However minute, this improvement is
perceptible as new grid points are added to the independent realizations.
Equally perceptible is the convergence property of the running average,
reaching faster an asymptotic value upon addition of grid points. Congruent
with what was established previously, the present set of results
 {emphasizes} that the SADM method cannot recover
the value  {derived} from the least-squares
equations. They further suggest that convergence towards the latter value will
only be achieved in the limit where the number of added points
 {coincides} with the total number of grid points
minus the number of parameters to be determined
--- {\em i.e.} one unique realization.

\begin{figure}[ht]
\begin{center}
\rotatebox{0}{\resizebox{!}{0.3\linewidth}{%
   \includegraphics{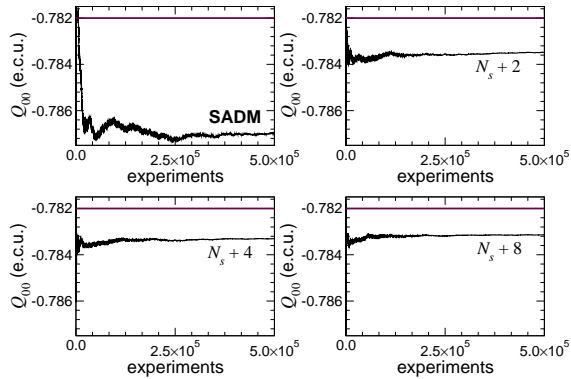}}}
\end{center}

\caption{Running average of the point charge, $Q_{00}$, borne by the oxygen
atom of water ($N_s$ = 2 parameters) as a function of the number of independent
realizations, wherein systems of 2 $\times$ 2 (SADM), 4 $\times$ 2, 6 $\times$
2 and 10 $\times$ 2 linear equations are solved. The thick, dark horizontal
line at $Q_{00}$ = $-$0.782 e.c.u. corresponds to the solution of the
least-squares problem. \label{averages} }
\end{figure}

Not too surprisingly, closer examination of the corresponding charge
distributions in Figure~{\ref{distributions}} reveals that as additional grid
points are added to the individual realizations, not only does the width of
these distributions narrow down, but the latter are progressively reshaped. As
was conjectured in \cite{chipot-98}, the SADM algorithm yields Cauchy
distributions, which is apparent from Figure~{\ref{distributions}}. Improvement
of the method alters the form of the probability function, now closer to a
normal distribution. Interestingly enough, the slightly skewed shape of the
distributions, particularly visible on their left-hand side --- as a probable
manifestation of the incompleteness of the electrostatic model, precludes
perfect enveloping by the model distributions, either Cauchy-- or
Gaussian--like.

\begin{figure}[ht]
\begin{center}
\rotatebox{0}{\resizebox{!}{0.3\linewidth}{%
   \includegraphics{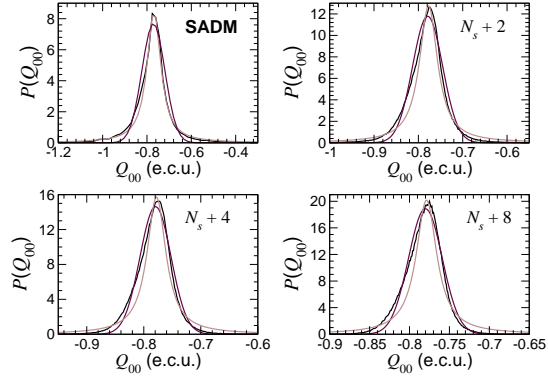}}}
\end{center}

\caption{Normalized distributions of the charge, $Q_{00}$, borne by the oxygen
atom of water ($N_s$ = 2 parameters) obtained from 500,000 independent
realizations, wherein systems of 2 $\times$ 2 (SADM), 4 $\times$ 2, 6 $\times$
2 and 10 $\times$ 2 linear equations are solved (black curves). The light and
dark curves correspond, respectively, to numerically fitted Cauchy and Gaussian
distributions. \label{distributions} }
\end{figure}

Put together, the present computations reinforce the conclusions drawn
hitherto, contradicting in particular the illegitimate assumption that the SADM
and the least-squares solutions might coincide \cite{chipot-98}. From a
numerical standpoint, however, the results obtained from both strategies appear
to be reasonably close, thereby warranting that the SADM algorithm should not
be obliterated, as it constitutes a valuable pedagogical tool for assessing the
appropriateness of electrostatic models.

\section{Conclusion}
 In this work, a probabilistic approach to high-dimensional least-squares
  approximations has been developed. Originally inspired by the SADM
  method introduced for the derivation of distributed atomic multipoles
  from the quantum-mechanical electrostatic potential, this novel approach
  can be generalized to a wide class of least-squares problems, yielding
  convergent and efficient numerical schemes in those cases where the space
  of approximation is very large or where the problem is ill-conditioned.

  This novel approach constitutes a marked improvement over the SADM method.
  Complete analysis of the numerical algorithm in general cases, in
  terms of both computational effort and optimal error estimation, relies
  on open and difficult issues prevalent to random matrix problems.


\appendix

\section{Proof of Lemma~\ref{Lborne}}
\label{sec:lemma}

We denote by $\Norm{\cdot}{F}$ the Schur-Frobenius norm on $n\times m$ matrices
{$$
\Norm{A}{F}^2=\sum_{i=1}^m\sum_{j=1}^n a_{ij}^2.
$$
where $A = (a_{ij})_{1\leq i \leq n,1\leq
  j\leq m}$.}
With this notation, we have
\begin{equation}
  \label{eq:ineg2}
  \Norm{(A^TA)^{-1}A^T}{F}^2
  =\mbox{Tr}((A^TA)^{-1}A^TA(A^TA)^{-1})=\mbox{Tr}((A^TA)^{-1})
  \leq \frac{n}{s_1(A)}.
\end{equation}
{In addition, for any} $v=(v_1,\ldots,v_m)\in\R^m$, we have
\begin{align}
  \Norm{Av}{2}^2 & =\sum_{i=1}^n\sum_{1\leq j,k\leq
    m}a_{ij}v_ja_{ik}v_k \notag \\ & \leq \frac{1}{2}\Big(
  \sum_{i,j,k}a_{ij}^2|v_jv_k|+\sum_{i,j,k}a_{ik}^2|v_jv_k|\Big)
  \notag \\ &
  \leq\Norm{A}{F}^2\Norm{v}{1}\Norm{v}{\infty}\leq\Norm{A}{F}^2\Norm{v}{1}^2,
  \label{eq:ineg1}
\end{align}
where we used the inequality $|b_{ij}b_{ik}| \leq \frac12 (b_{ij}^2 + b_{ik}^2)$.

{With the notation of Lemma~\ref{Lborne},} using \eqref{eq:ineg1} and \eqref{eq:ineg2}, we have
$$
\Norm{R(X)G(X)}{2} \leq \sqrt{n}\, s_1(\Gamma(X))^{-1/2} \Norm{G(X)}{1}.
$$
Taking the expectation and using H\"older's inequality, we get
\begin{equation}
\label{eq:pf-lemma}
\E\Norm{R(X)G(X)}{2} \leq \sqrt{n}\sqrt{K_q(\Gamma)} \Big(\E \Norm{G(X)}{1}^p\Big)^{1/p}.
\end{equation}
Now, $Y\mapsto(\E |Y|^{p})^{1/p}$ defines a norm on the set of random vectors
on $\Omega$ with finite $p$-th-order moment. We,
 {hence}, obtain
\begin{align*}
  \big(\E\Norm{G(X)}{1}^p\big)^{1/p} & =\Big(\E\Big(\sum_{i=1}^m
  |g(X^{(i)})|\Big)^{p}\Big)^{1/p} \\
  & \leq\sum_{i=1}^m\big(\E
  |g(X^{(i)})|^p\big)^{1/p}=m\Norm{g}{\Lr^p(\Omega)},
\end{align*}
and this yields {Lemma~\ref{Lborne}~(a).

Lemma~\ref{Lborne}~(b)} is obtained from similar computations.

\section{Estimates on random matrices}

{\subsection{On the finiteness of $K_q(\Gamma)$}}
\label{sec:moment}

{As seen in Propositions~\ref{Tmoments} and~\ref{prop:comp-lsquares}, the convergence {and the consistency of
  our algorithm} rely on assumptions of the form $K_q(\Gamma)<+\infty$, where $K_q(\Gamma)$ is given
by~(\ref{EKmq}). These assumptions correspond to {the finiteness of a negative moment of the random variable
  $s_1(\Gamma(X))$. The following proposition {provides} a condition on the distribution of $s_1(\Gamma(X))$
  to ensure such integrability properties.}}

\begin{proposition}
  \label{prop:mom}
  Let $Y$ be a random variable satisfying the following estimate:
 There exist constants $\delta > 0$ and $\gamma > 0$ such
  that
  \begin{equation}
    \label{Econdsn}
    \forall\epsilon \geq 0,\quad \P(Y \leq \epsilon)
    \leq (\delta \epsilon)^{\gamma}.
  \end{equation}
  Then, for any $0 < r < \gamma$,  \begin{equation}
    \label{eq:mom_sigma}
    \E(Y^{-r})\leq
    \frac{\delta^r}{1-r/\gamma}.
  \end{equation}
\end{proposition}

\begin{Proof}
{The proof is based on} the following integration by parts, where
$\int_0^\infty h(x)\dd \P(Y\in[0,x))$ denotes the
Stieltjes integral of the measurable function $h$ with respect to
the Stieltjes measure on $[0,\infty)$ associated with the
non-decreasing function $x\mapsto \P(Y\in[0,x))$.
\begin{align*}
  \E(Y^{-r}) & =\int_0^\infty
  x^{-r} \dd
  \P(Y\in[0,x)) \\[1ex]
  & =\int_0^\infty r x^{-r-1}
  \P(Y\in[0,x))
  \dd x \\[1ex]
  & \leq r
  \int_0^\infty x^{-r-1}((\delta x)^\gamma\wedge 1) \dd x.
\end{align*}
If $r < \gamma$,
\begin{equation*}
  \int_0^\infty x^{-r-1}((\delta x)^\gamma\wedge 1) \dd x=
  \frac{1}{r}
  \Big(\frac{\delta^r}{1-r/\gamma}
  \Big),
\end{equation*}
which entails~(\ref{eq:mom_sigma}).
\end{Proof}

{Of course, the property~\eqref{Econdsn} can be strongly problem-dependent.} In general situations, this is
related to difficult problems on random matrices, which, to our knowledge, have not been solved
yet. {However, explicit computations are possible in the specific instance where $\Gamma(X)^T\Gamma(X)$ has
the Wishart distribution {(see Subsections~\ref{sec:Wishart} and below).}
\medskip

{In the case where the matrix $\Gamma^T(X)\Gamma(X)$ is not a.s.\ invertible, the method described in
  Section~\ref{sec:sigma} consists in taking expectations conditionally to $\{s_1(\Gamma(X))>\sigma\}$ for
  some $\sigma>0$. A quantitative analysis of our method relies on estimates on $K_q^\sigma(\Gamma)$ defined
  in~\eqref{eq:def-K-sigma} (see Propositions~\ref{prop:consist2} and~\ref{Tmoments2}).} The following result
generalizes Proposition~\ref{prop:mom} to the case where $\sigma > 0$. Its proof is very similar to {that} of
Proposition~\ref{prop:mom}. We {will, hence,} omit it {here}.

\begin{proposition}
  \label{prop:mom2}
  Fix $\sigma\geq 0$ and assume that random variable $Y$ satisfies the
  following estimate: There exist constants $\delta$ and $\gamma$ such
  that
  \begin{equation}
    \label{Econdsn2}
    \forall\, \epsilon\geq\sigma,\quad \P(Y \leq \epsilon)
    \leq (\delta \epsilon)^{\gamma}.
  \end{equation}
  Then, for any $r\not=\gamma$, if $0<\sigma<\delta^{-1}$,
  \begin{equation}
    \label{eq:mom_sigma2}
    \E(Y^{-r}\mid Y\geq\sigma)\leq
    \frac{\delta^r}{1-(\delta\sigma)^\gamma}\Big(
    \frac{1}{1-r/\gamma}+\frac{(\delta\sigma)^{-(r-\gamma)}}{1-\gamma/r}\Big).
  \end{equation}
\end{proposition}

{This result is used to obtained explicit estimates on our algorithm in the case where the matrix
  $\Gamma(X)$ has sub-Gaussian entries (see Sections~\ref{sec:subgaussian} and~\ref{sec:sub-gaussian-pf}).}

{\subsection{Explicit estimates in the Wishart case}}
\label{sec:Wishart-pf}

The goal of this section is to prove the following result.

\begin{proposition}
  \label{prop:Wishart-algo} In {the Wishart case {(see Section~\textup{\ref{sec:Wishart}}),} assume that
  $f\in\Lr^p(\Omega)$ and $\rho(\bar{\beta})\in \Lr^p(\Omega)$ with $p > 2$. Then, the convergence of the
  algorithm (in the sense that $\textup{Tr}(\textup{Cov}(\beta))<+\infty$, see
  Proposition~\textup{\ref{Tmoments}}) holds
  if
  $$
  m > n + \frac{p+2}{p-2}.
  $$
  In the case where $\rho(\bar{\beta})\in \Lr^\infty(\Omega)$, this condition corresponds to $m\geq n+2$ and
  the computational cost of the algorithm is bounded by
  $$
  C\varepsilon^{-2} \, \frac{ n^3m^4\log n }{(m-n+1)(m-n-1)} \Norm{\rho(\bar\beta)}{\Lr^\infty}^2,
  $$
  where $\varepsilon$ is the required precision and the constant $C$ is independent of $n$ and $m$. The
  optimal value $m^*$ of $m$ in the previous bound satisfies
  $$
  m^*\sim 2n
  $$
  when $n\rightarrow+\infty$, and the corresponding computational cost is bounded by
  $$
  C'\varepsilon^{-2} n^5\log n  \Norm{\rho(\bar\beta)}{\Lr^\infty}^2.
  $$
  In addition, the consistency error of Proposition~\ref{prop:comp-lsquares} is bounded by
  $$
  C''n\Norm{\rho(a)}{\Lr^\infty}.
  $$}
\end{proposition}

{Our computations are based on the following estimate on the law of the smallest eigenvalue of Wishart
  matrices~\cite[Lemma~3.3]{chen-dongarra-05}, which reads with our notation as follows.} For all $m\geq n\geq
2$, let 
$$
k=m-n+1.
$$
The density $p(x)$ of $s_1(\Gamma(X))$ {then} satisfies
\begin{equation}
  \label{eq:densite-s_1}
  L_{n,m}e^{-nx/2}x^{\frac{k}{2}-1}\leq p(x)\leq
  L_{n,m}e^{-x/2}x^{\frac{k}{2}-1},\quad\forall x>0,
\end{equation}
where
\begin{equation}
  \label{eq:def-lnm}
  L_{n,m}=\frac{2^{\frac{k}{2}-1}\Phi(\frac{m+1}{2})}
  {\Phi(\frac{n}{2})\Phi(k)},
\end{equation}
where $\Phi$ is the Gamma function, defined for all $x>0$
by
$$
\Phi(x)=\int_0^{+\infty}e^{-t}t^{x-1}\dd t.
$$

\begin{lemma}
  \label{lem:Wishart}
   For all $m\geq n\geq 2$, the random variable $Y=s_1(\Gamma(X))$
   satisfies~\eqref{Econdsn} for
   \begin{gather*}
     \gamma=\frac{m-n+1}{2}=\frac{k}{2} \\
     \mbox{and}\quad\delta=e^2\frac{m}{k^2}.
   \end{gather*}
   Moreover, the constant $\gamma$ above is the smallest such
   that~\eqref{Econdsn} holds for all $\varepsilon>0$ for some
   constant $\delta$.
\end{lemma}

\begin{Proof}
  {This result is based on} the following bounds for the
  Gamma function~\cite[Lemma 2.7]{chen-dongarra-05}. For all $x>0$,
  $$
  \sqrt{2\pi}\:x^{x+\frac{1}{2}}e^{-x}<\Phi(x+1)=x\Phi(x)
  <\sqrt{2\pi}\:x^{x+\frac{1}{2}}e^{-x+\frac{1}{12x}}.
  $$
  These inequalities can be plugged into~(\ref{eq:def-lnm}) to get
  that, for all $\varepsilon>0$,
  \begin{align*}
    \P(s_1(\Gamma(X))\leq\varepsilon) & \leq
    \frac{n\:2^{\gamma-1}\Phi(\frac{m-1}{2}+1)}
    {\Phi(\frac{n}{2}+1)2\gamma\Phi(2\gamma)}
    \:\gamma\int_0^{\varepsilon}x^{\gamma-1}\dd x \\
    & \leq
    \frac{\frac{n}{2}2^\gamma\sqrt{2\pi}(\frac{m-1}{2})^{\frac{m}{2}}
      e^{-\frac{m-1}{2}+\frac{1}{6(m-1)}}}
    {\sqrt{2\pi}(\frac{n}{2})^{\frac{n+1}{2}}e^{-\frac{n}{2}}
      \sqrt{2\pi}(2\gamma)^{2\gamma+\frac{1}{2}}e^{-2\gamma}}\:
    \varepsilon^{\gamma} \\
    & \leq \frac{e^{1+\frac{1}{6(m-1)}}}{2\sqrt{\pi\gamma}}\:
    \frac{(\frac{m-1}{2})^{\frac{m}{2}}}{(\frac{n}{2})^{\frac{n-1}{2}}}\:
    \Big(\frac{2e^2\varepsilon}{4\gamma^2e}\Big)^{\gamma}.
  \end{align*}
  Now,
  \begin{align*}
    \frac{(\frac{m-1}{2})^{\frac{m}{2}}}{(\frac{n}{2})^{\frac{n-1}{2}}}
    =\Big(\frac{m-1}{n}\Big)^{\frac{n-1}{2}}\Big(\frac{m-1}{2}\Big)^{\gamma} &
    \leq\Big(1+\frac{k-2}{n}\Big)^{\frac{n}{2}}\Big(\frac{m-1}{2}\Big)^{\gamma} \\ &
    \leq e^{\frac{k-2}{2}}\Big(\frac{m-1}{2}\Big)^{\gamma}
    =\frac{1}{e}\Big(\frac{e(m-1)}{2}\Big)^\gamma.
  \end{align*}
  Combining this inequality with the facts that $m-1\geq 1$ and
  $\gamma\geq 1/2$ yields
  $$
  \P(s_1(\Gamma(X))\leq\varepsilon)
  \leq\frac{e^{1/6}}{\sqrt{2\pi}}\Big(\frac{e^2(m-1)}{4\gamma^2}\:
  \varepsilon\Big)^\gamma\leq\Big(\frac{e^2m}{k^2}\:\varepsilon\Big)^\gamma.
  $$

  Because of~(\ref{eq:densite-s_1}), we have that $p(x)\sim
  L_{n,m}x^{\gamma-1}$ as $x\rightarrow 0$. Therefore, one easily sees
  that $\gamma=k/2$ is the minimal value of $\gamma$
  for~(\ref{Econdsn}) to holds.
\end{Proof}

Using this result and Proposition~\ref{prop:mom}, we immediately
obtain the following:
\begin{proposition}
Let $m > n$ be given and assume that the random matrix $\Gamma(X)^T\Gamma(X)$
associated with the function $\Gamma$ defined in \eqref{E4} follows a Wishart
distribution. Let $q$ be such that
\begin{equation}
\label{Econdq}
1 \leq q < k = m - n +1.
\end{equation}
Then we have
\begin{equation}
\label{EmomWishart}
K_q(\Gamma) \leq \frac{e^2m}{k^2}\Big(1 - \frac{q}k\Big)^{-\frac{2}{q}}
\end{equation}
where $K_q$ is defined in \eqref{EKmq}.
\end{proposition}

Combining this result and the result of Proposition \ref{Tmoments}, if
$\rho(\bar\beta) \in \Lr^p(\Omega)$ with $p > 2$ , the convergence of the
algorithm is ensured if $K_{q} < \infty$ in \eqref{eq:var*1} with
$2p^{-1}+2q^{-1} = 1$. This means, (see \eqref{Econdq})
$$
2 \leq \frac{2p}{p-2} < m - n +1
$$
or equivalently
$$
m > n + \frac{p+2}{p-2}.
$$

Assume still that $\rho(\bar\beta) \in \Lr^\infty(\Omega)$. Using
\eqref{EmomWishart} with $q = 2$,  {it can be
seen} in view of \eqref{EKost} that the cost of the algorithm is bounded by
$$
C\varepsilon^{-2}  n^3m^3\log n \, \frac{ m }{k^2 }\Big(1 - \frac{2}k\Big)^{-1}
\Norm{\rho(\bar\beta)}{\Lr^\infty}^2.
$$
for some constant $C$ independent of $n$ and $m$.
Using the notation $\gamma = k/2$, we  can rewrite this cost in term of $\gamma$ as
$$
C'\varepsilon^{-2}n^3\log
  n \, \frac{(n+2\gamma-1)^4}{\gamma^(\gamma-1)}\Norm{\rho(\bar\beta)}{\Lr^\infty}^2.
$$
To determine the optimal choice of $m$,  {let us}
now try to find the optimal number $\gamma$ that
 {minimizes} this cost. The derivative of this
expression with respect to $\gamma$ has the same sign as
$$
  8\gamma(\gamma-1)-(n+2\gamma-1)(2\gamma-1)
  =4\gamma^2-2(n+2)\gamma+n-1.
$$
Since this quantity is negative if $\gamma=1/2$, the only root of this
polynomial  {greater} than $1$ is given by
$$
\gamma^*=\frac{n+2+\sqrt{n^2 + 8}}{4},
$$
which is the optimal choice of $\gamma$ in terms of computational
 {effort}. This
 {yields} an optimal choice $m^* \sim 2 \gamma^*
+ n - 1$. Note that for large $n$, we have $\gamma^* \sim n/2$ and $m^* \sim 2
n$.

With this optimal choice, the  {computational}
cost of the algorithm can be written as
\begin{equation*}
  C_n \varepsilon^{-2}\,\Norm{\rho(\bar\beta)}{\Lr^\infty}^2\quad \mbox{with}\quad C_n \sim  C n^5\log n
  \quad\mbox{as}\quad n \to +\infty. 
\end{equation*}
\medskip

{Considering a} similar calculation with $q=1$,
we  {can} easily see that the consistency error
of Proposition~\ref{prop:consist} for this choice of parameters can be bounded
by
$$
C'_n \Norm{\rho(a)}{\Lr^\infty}\quad\mbox{with}\quad
C'_n \sim C' n \quad\mbox{as}\quad n \to +\infty.
$$

{\subsection{Explicit estimates in the  sub-Gaussian case}}
\label{sec:sub-gaussian-pf}

{The goal of this section is to prove the following result.}

{\begin{proposition}
    \label{prop:subgaussian} In the sub-Gaussian case {(see Section~\textup{\ref{sec:subgaussian}})}, assume
    that $f\in\Lr^\infty(\Omega)$ and $\rho(\bar{\beta})\in \Lr^\infty(\Omega)$. Then, there exists explicit
    constants $A$ and $B$ such that, if
  $$
  \sigma=\frac{B^2m^2(\sqrt{m}-\sqrt{n-1})^2}{A(m-n+1)^2}e^{-2Bm/(m-n+1)},
  $$
  the computational cost of our algorithm is bounded by
  \begin{equation}
    \label{eq:cost-subg}
    C\varepsilon^{-2} \, \frac{ n^3m^4\log n }{(m-n+1)(m-n-1)} \Norm{\rho(\bar\beta)}{\Lr^\infty}^2,    
  \end{equation}
  where $\varepsilon$ is the required precision and the constant $C$ is independent of $n$ and $m$. Again, the
  optimal value $m^*$ of $m$ in the previous bound satisfies $m^*\sim 2n$ as $n\rightarrow+\infty$.\\
  For such a choice of $m$, we obtain
  \begin{equation}
    \label{eq:comp-sigma}
    \sigma\sim C'\,n    
  \end{equation}
  for an explicit constant $C'$.
\end{proposition}}

{With our notations, Theorem~1.1 of~\cite{rudelson-vershynin-08} writes as follows: there exist explicit}
constants $A$ and $B$ depending only on $R$ such that, for all $m\geq n$ and all $\epsilon>0$,
\begin{equation}
  \label{eq:R-V}
  \P\Big(s_1(\Gamma(X))\leq\epsilon(\sqrt{m}-\sqrt{n-1})^2\Big)
  \leq (A\epsilon)^{(m-n+1)/2}+e^{-Bm}.
\end{equation}
 {Writing just like} in
Subsection~\ref{sec:Wishart-pf} $k$ for $m-n+1$,
 {it can be seen} that
$$
\begin{array}{rcl}
  \P(s_1(\Gamma(X))\leq\varepsilon)
  &\leq & \displaystyle \Big(\frac{\sqrt{A\varepsilon}}{\sqrt{m}-\sqrt{n-1}}\Big)^{k}
  +\big(e^{-Bm/k}\big)^k \\[2ex]
  \
  &\leq & \displaystyle
  \Big(\frac{\sqrt{A\varepsilon}}{\sqrt{m}-\sqrt{n-1}} +e^{-Bm/k}\Big)^k.
  \end{array}
$$
Eq.~(\ref{Econdsn2}), therefore, holds for $Y=s_1(\Gamma(X))$ and
\begin{gather*}
  \sigma
  \geq\sigma_0:=\frac{B^2m^2(\sqrt{m}-\sqrt{n-1})^2}{k^2A}e^{-2Bm/k},
  \\ \delta=\frac{(1+k/Bm)^2A}{(\sqrt{m}-\sqrt{n-1})^2} \\
  \mbox{and}\quad \gamma=\frac{k}{2}.
\end{gather*}
Note that, since $\delta\sigma_0=(1+Bm/k)^2e^{-2Bm/k}< 1$, the inequality in
(\ref{Econdsn2}) is not trivial and  {supplies}
some information on the law of $s_1(\Gamma(X))$.

As in Subsection~\ref{sec:Wishart-pf}, the inequality (\ref{eq:mom_sigma2}) can
 {be combined} with the results
 {of} Propositions~\ref{Tmoments2} and~\ref{prop:consist2} to obtain  {a}
precise error  {estimate} and convergence bounds
in this case.

{Such} computations are, {however}, {cumbersome} because the optimal choice of $\sigma$ cannot be {determined
  explicitly}.
{Taking $\sigma=\sigma_0$ as in Proposition~\ref{prop:subgaussian} and} under the assumption that
$\rho(\bar{\beta}^\sigma)\in\Lr^\infty(\Omega)$, because of Proposition~\ref{prop:mom2}, the computational
cost is smaller than
$$
C\varepsilon^{-2}n^3\log n\Norm{\rho(\bar{\beta}^\sigma)}{\Lr^\infty}^2
\frac{m^3\delta}{(1-(\delta\sigma)^\gamma)^2}
\Big(\frac{1}{1-1/\gamma}+\frac{(\delta\sigma)^{\gamma-1}}{1-\gamma}\Big)
$$
If one assumes that
$(\sigma_0\delta)^\gamma\rightarrow 0$ as $n\rightarrow+\infty$,
observing that
$$
\delta=\frac{(\sqrt{m}+\sqrt{n-1})^2(1+k/Bm)^2A}{k^2}\leq
\frac{Cm}{k^2},
$$
the cost is bounded from above by
$$
C\varepsilon^{-2}n^3\log n
\frac{m^4}{\gamma(\gamma-1)}\Norm{\rho(\bar{\beta}^\sigma)}{\Lr^\infty}^2
$$
for $\gamma>1$. We recognize the same cost as in Subsection~\ref{sec:Wishart}.
The optimal choice of $\gamma$,  {therefore},
behaves as $n/2$ as $n\rightarrow +\infty$ --- and for this choice we indeed
have $(\sigma_0\delta)^\gamma\rightarrow 0$, which validates the
previous computation. Therefore, for this choice of parameters,
 {the} cost is bounded by
$$
C\varepsilon^{-2}n^5\log n\Norm{\rho(\bar{\beta}^\sigma)}{\Lr^\infty}^2
$$
for some constant $C>0$.

One can check that any other choice of $\sigma$ {yields} the same order in $n$ as $n\rightarrow+\infty$,
{should one choose} $\gamma\sim n/2$.  \medskip

{It ought to be noted} that these bounds do not allow one to {pick} $\sigma=0$. As far as we know, {this}
seems to be an open and difficult question to prove that~(\ref{eq:R-V}) holds without the {right-hand-side},
additive term $e^{-Bm}$. In particular, it requires additional assumptions to hold --- {{\it e.g.}} random
variable $h(Y)$, where $Y$ has law $\nu$, has no atom, \emph{i.e.}\ that $\nu(\{h=y\})=0$ for all $y\in\R$
(otherwise, the matrix $\Gamma(X)$ could have $m-n$ identical rows, and, {thus}, have a rank less than $n$,
with non-zero probability).

\bigskip

\bibliographystyle{plain}

\end{document}